\magnification = 1200
\showboxdepth=0 \showboxbreadth=0

\baselineskip 14pt
\parskip3pt
\def\qed{\hfill\vrule height6pt width6pt depth0pt}

\def\ss{\smallskip}
\def\ms{\medskip}
\def\bs{\bigskip}

\def\cl{\centerline}

\def\nind{\noindent}

\def\ref#1#2{\nind\hangindent.5in\hbox to .5in{#1\hfill}#2}
\def\reff#1#2{\nind\hangindent.8in\hbox to .8in{\bf #1\hfill}#2\par}
\def\refd#1#2{\nind\hangindent.8in\hbox to .8in{\bf #1\hfill {\rm
--}}#2\par}
\def\pmb#1{\setbox0=\hbox{#1}
\kern-0.025em\copy0\kern-\wd0
\kern.05em\copy0\kern-\wd0
\kern-.025em\raise.0433em\box0}
\def\ca#1{{\cal #1}}

\def\Cal{\ca}

\def\frac#1#2{{#1\over#2}}

\def\text#1{\rm#1}

\outer\def\stmnt  #1. #2\par{\medbreak
\noindent{\bf#1.\enspace}{\sl#2}\par
\ifdim\lastskip<\medskipamount \removelastskip\penalty55\medskip\fi}

\def\newline{\hfill\break}
\def\:{\,:\,}

\def\({\left(}                   
\def\){\right)}

\def\[{\left[}                   
\def\]{\right]}

\def\lan{\langle}
\def\ran{\rangle}

\def\ci{\subset}

\def\fy{\infty}

\def\del{\partial}

\def\Om{\Omega}

\cl{\bf The Aronsson equation for absolute minimizers of ${L^\fy}$-functionals}
\cl{\bf associated with vector fields satisfying H\"ormander's condition}
\bs
\cl{Changyou Wang}
\ss
\cl{Department of Mathematics, \ University of Kentucky}
\cl{Lexington, KY 40506}  
\bs
\nind{\bf Abstract}. {\it Given a Carnot-Carath\'eodory metric
space $(R^n, d_{\hbox{cc}})$ generated by vector fields $\{X_i\}_{i=1}^m$ 
satisfying H\"ormander's condition, we prove in theorem A that any absolute minimizer 
$u\in W^{1,\fy}_{\hbox{cc}}(\Om)$ to $F(v,\Om)=\sup_{x\in\Om}f(x,Xv(x))$ 
is a viscosity solution to the Aronsson equation (1.6), under suitable conditions on $f$. 
In particular, any AMLE is a viscosity solution to the subelliptic $\fy$-Laplacian equation (1.7). 
If the Carnot-Carath\'edory space is a Carnot group ${\bf G}$
and $f$ is independent of $x$-variable, we establish in theorem C the uniquness of viscosity solutions
to the Aronsson equation (1.13) under suitable conditions on $f$. 
As a consequence, the uniqueness of both AMLE and viscosity solutions to the subelliptic
$\fy$-Laplacian equation is established in ${\bf G}$.}
\bs
\nind{\S1}. Introduction
\ss
Variational problems in $L^\fy$ are very important because of both
its analytic difficulties and their frequent appearance in applications,
see the survey article [B] by Barron . 
The study began with Aronsson's papers [A1, 2]. 
The simplest model is to consider minimal Lipschitz extensions (or MLE):
for a bounded, Lipschitz domain $\Om\ci R^n$ and $g\in \hbox{Lip}(\Om)$, find  
$u\in W^{1,\fy}(\Om)$, with $u|_{\del\Om}=g$, such that
$$\|Du\|_{L^\fy(\Om)}\le \|Dw\|_{L^\fy(\Om)},
\ \forall w\in W^{1,\fy}_0(\Om), \hbox{ with } w|_{\del\Om}=g. \eqno(1.1)$$  
Since MLE's may be neither unique nor smooth, Aronsson [A1] introduced the
notation of absolutely minimizing Lipschitz extensions(or AMLE for short),
and proved that any $C^2$ AMLE solves the $\fy$-Laplacian equation
$$\Delta_\fy u: =-\sum_{i j =1}^n {\del u\over\del x_i}{\del u\over\del x_j}
{\del^2 u\over\del x_i\del x_j} = 0,  \ \ \hbox{ in }\ \Om. \eqno(1.2)$$
However, (1.2) is a highly nonlinear and highly degenerate PDE and
may not have $C^2$ solutions in general. This issue was finally settled by
Jensen [J2], who not only established the equivalence between the AMLE property and  
the solution to eqn.(1.2) in the viscosity sense, which was first introduced
by Crandall-Lions [CL](see also Crandall-Ishii-Lions [CIL]), but also proved
the uniqueness of viscosity solutions to eqn.(1.2) with the Dirichlet boundary
value. The remarkable analysis of [J2] involves approximation by $p$-Laplacian
and Jensen's earlier work [J1] on the maximum principle for semiconvex functions. 
The reader can consult with Evans [E] and  Lindqvist-Manfredi [LM] for qualitative
estimates on $\fy$-harmonic functions.
Crandall-Evans-Gariepy [CEG] developped the comparison principle
of the eqn.(1.2) with {\it cones}, which are solutions of the eiknonal
equation of forms $a+b|x-x_0|$, and gave an alternative, direct proof 
of the equivalence between AMLE and viscosity solution to eqn.(1.2).
Furthermore, Crandall-Evans [CE] has utilized this property in their
study on the regualarity issue of $\fy$-harmonic functions. 
Recently, Barron-Jensen-Wang [BJW] considered general $L^\fy$-functionals
$$F(u,\Om):=\sup_{x\in\Om}f(x,u(x),Du(x)), \ \forall u \in W^{1,\fy}(\Om).$$
and proved, under suitable conditions, that any absolute minimizer
of $F(\cdot,\Om)$ is a viscosity solution to the Aronsson-Euler equation
$$-\sum_{i=1}^n f_{p_i}(x,u(x),Du(x)) {\del \over\del x_i}
(f(x,u(x),Du(x))) = 0,  \ \ \hbox{ in }\ \ \Om. 
\eqno(1.3)$$  
Shortly after [BJW], Crandall [C] was able to give an elegant proof of an improved
version of [BJW]. Through [BJW] [C], it becomes more clear that the classical
solution to the Hamilton-Jacobi equation $f(x,\phi(x), D\phi(x))-c = 0$
plays important roles in this analysis. 
 
Since the notion of AMLE can easily be formulated in any metric space, it is
a very natural and interesting problem to study AMLE in spaces with 
Carnot-Carath\'eodory metrics, which include Riemannian manifolds and 
Subriemannian manifolds (e.g., Heisenberg groups, Carnot groups, or more generally
H\"ormander vector fields, etc). There have been several works
done in this direction. For example, Juutinen [J] extended the main theorems of
[J2] into Riemannian manifolds. Bieske [B1,2] was able to prove that,
on the Heisenberg group ${\bf H}^n$ or a Grushin type space, an AMLE is equivalent 
to a viscosity solution to the subelliptic $\fy$-Laplacian equation, and the
uniqueness of both AMLE and viscosity solution to the subelliptic $\fy$-Laplacian 
equation. Inspired by Crandall's argument [C], Bieske-Capogna, in a recent 
preprint [BC], proved that any AMLE is a viscosity solution to the subelliptic 
$\fy$-Laplacian equation for any Carnot group, where the conclusion was also proved for any AMLE,
which is horizontally $C^1$, corresponding to those Carnot-Carath\'edory metrics 
associated to {\it free} systems of H\"ormander's vector fields.  

In this paper, we are mainly interested in the derivation of Euler equation
of AMLE and its uniqueness issue for any Carnot-Carath\'edory metric space
generated by vector fields satisfying H\"ormander's condition. In this direction,
we are able to prove that any AMLE is a viscosity solution to the subelliptic
$\fy$-Laplacian equation. Moreover, if the vector fields are horizontal 
vector fields associated with a Carnot group, then we establish the uniqueness
for both AMLE and viscosity solution to the Euler equation. In fact, these
conclusions are consequences of general theorems A and C below. 

In order to state our results, we first recall some preliminary facts.
\ss
\nind{\bf Definition 1.1}. For a bounded domain $\Om\ci R^n$ and $m\ge 1$,
$\{X_i\}_{i=1}^m \ci C^2(\Om, R^n)$ are vector fields satisfying H\"ormander's
condition, if there is a step $r\ge 1$ such that, at any $x\in \Om$, 
$\{X_i\}_{i=1}^m$ and all their commutators up to at most order $r$ generate $R^n$.  

Now we recall from [NSW] the Carnot-Carath\'edory distance, denoted as 
$d_{\hbox{cc}}$, generated by $\{X_i\}_{i=1}^m$ on $\Om$:  for any $p,q\in \Om$, 
$$d_{\hbox{cc}}(p,q)=\inf_{A(\delta)}\delta, \eqno(1.4)$$ 
where $$A(\delta):=\{r:[0,\delta]\to \Om \ |\ r(0)=p,\ r(\delta)=q,
r'(t)=\sum_{i=1}^m a_i(t)X_i(r(t)) \hbox{ with }\sum_{i=1}^m a_i^2(t)\le 1\}.$$ 
Moreover, $d_{\hbox{cc}}$ satisfies: for each compact set $K\ci\ci \Om$,
$$C_K^{-1}\|x-y\|\le d_{cc}(x,y)\le C_K\|x-y\|^{1\over r}, \ \forall x, y\in K, $$
where $\|\cdot\|$ is the Euclidean norm on $R^n$.
For $u:\Om\to R$, denote $Xu:=(X_1 u,\cdots, X_m u)$ as
the horizontal gradient of $u$. 
For $1\le p\le\fy$, the horizontal Sobolev space is defined by
$$W_{\hbox{cc}}^{1,p}(\Om)
:=\{u: \Om \to R \  |\ \  \|u\|_{W^{1,p}_{\hbox{cc}}(\Om)}
\equiv\|u\|_{L^p(\Om)}+\|Xu\|_{L^p(\Om)}<\fy 
\}.$$
The Lipschitz space, with respect to the metric $d_{\hbox{cc}}$, is defined by
$$\hbox{Lip}^{\hbox{cc}}(\Om):=\{u: \Om\to R \  |
\ \|u\|_{{\hbox {Lip}}^{\hbox{cc}}(\Om)}\equiv\sup_{x,y\in \Om, x\not= y}
{|u(x)-u(y)|\over d_{\hbox{cc}}(x,y)}<\fy\}.$$
It was proved by [GN] (see also [FSS]) that
$u\in\hbox{Lip}^{\hbox{cc}}(\Om)$ iff $u\in W^{1,\fy}_{\hbox{cc}}(\Om)$.

Now we recall the definition of absolute minimizers of $L^\fy$-functionals
over $W^{1,\fy}_{\hbox{cc}}(\Om)$.
\ss
\nind{\bf Definition 1.2}. For any integrand function $f:\Om\times R^m\to R_+$,
let
$$F(v,\Om)=\sup_{x\in\Om}f(x,Xu(x)), \ \forall v\in W^{1,\fy}_{\hbox{cc}}(\Om).$$
A function $u\in W^{1,\fy}_{\hbox{cc}}(\Om)$ is an
{\it absolute minimizer} of $F(\cdot,\Om)$, if for any open subset $\tilde{\Om}\ci\Om$ and
$w\in W^{1,\fy}_{\hbox{cc}}({\tilde\Om})$, with $w=u$ on $\del {\tilde\Om}$, we have
$$F(u, {\tilde\Om})\le F(w, {\tilde\Om}). \eqno(1.5)$$
$u$ is called an {\it absolutely minimizing Lipschitz extension} (or AMLE), with respect
to the Carnot-Carath\'edory metric $d_{\hbox{cc}}$, if $u$ is an absolute minimizer
of $F(\cdot,\Om)$, with $f(x,p)=\sum_{i=1}^m p_i^2$ for $(x,p)\in\Om\times R^m$.

Formal calculations yield that an absolute minimizer $u\in W^{1,\fy}_{\hbox{cc}}(\Om)$
to $F(\cdot,\Om)$ satisfies the subelliptic Aronsson-Euler equation
$$-\sum_{i=1}^m X_i(f(x,Xu(x))) f_{p_i}(x,Xu(x)) = 0, \ \hbox{ in }\ \Om.
\eqno(1.6)$$
In particular, the Aronsson-Euler equation of an AMLE is the
subelliptic $\fy$-Laplacian equation
$$\Delta_\fy^{(X)} u: =-\sum_{i,j=1}^m X_i u X_j u X_i X_j u = 0,
\ \hbox{ in } \ \Om. \eqno(1.7)$$

In order to interpret an absolute minimizer (or AMLE respectively)
as a solution to the eqn. (1.6) (or (1.7) respectively), we recall the
concept of viscosity solutions by Crandall-Lions [CL] (see also [CIL]) 
of second order degenerate subelliptic PDEs.  

Let ${\Cal S}^m$ denote the set of symmetric $m\times m$ matrices, equipped
with the usual order. A function $A\in C(R^n\times R^m\times {\Cal S}^m)$ is
called degenerate subelliptic, if, for any $(x,p)\in R^n\times R^m$
$$A(x,p, M)\le A(x,s, N),  \ \forall M, N\in {\Cal S}^m , \ \hbox{ with } N\le M. \eqno(1.8)$$
Let $(D^2 u)^*\in {\Cal S}^m$ denote the {\it horizontal} hessian of $u$, defined by
$$(D^2u)^*_{ij}={1\over 2}(X_iX_j+X_jX_i)u,  \ \forall 1\le i, \ j \le m.$$
Now we have
\ss
\nind{\bf Definition 1.3}. For a degenerate subelliptic equation
$$A(x,Xu(x),(D^2u)^*(x)) = 0, \  \hbox{ in }\ \Om . \eqno(1.9)$$
A function $u\in C(\Om)$ is called a viscosity subsolution to eqn.(1.9),
if for any pair $(x_0,\phi)\in\Om\times C^2(\Om)$ such that
$x_0$ is a local maximum point of $(u-\phi)$ then we have
$$A(x_0,X\phi(x_0), (D^2\phi)^*(x_0))\le 0. \eqno(1.10)$$
A function $u\in C(\Om)$ is called a viscosity supersolution to eqn.(1.9)
if $-u$ is a viscosity subsolution to eqn.(1.9). Finally, a function 
$u\in C(\Om)$ is a viscosity solution to eqn.(1.9) if it is
both a viscosity subsolution and a viscosity supersolution to eqn.(1.9).

It is easy to check that both eqn. (1.6) and (1.7) are degenerate subelliptic.
Now we are ready to state our first theorem.
\ss
\nind{\bf Theorem A}. {\it Suppose that $u\in W^{1,\fy}_{\hbox{cc}}(\Om)$ 
is an absolute minimizer of $$F(v,\Om)=\sup_{x\in\Om}f(x,Xv(x)),$$ where   
$f\in C^2(\Om\times R^m,R_+)$ satisfies

\nind{(f1)} $f$ is quasiconvex in its second variable, i.e. for any $x\in\Om$,
$$f(x, tp_1+(1-t)p_2)\le\max\{f(x,p_1), f(x, p_2)\}, \ \forall p_1 \ p_2\in R^m, \ 0\le t\le 1.
\eqno(1.11)$$
\nind{(f2)} $f$ is homogeneous of degree $\alpha\ge 1$ and $f_p(0,0)=0$.

\nind  Then $u$ is a viscosity solution to the Aronsson-Euler equation
$$-\sum_{i=1}^m X_i(f(x,Xu(x)))f_{p_i}(x,Xu(x)) = 0, \  \ \hbox{ in }\ \Om. \eqno(1.12)$$}   

The ideas to prove theorem A are based on: (1) the observation of rewrite eqn.(1.12) into 
an euclidean form, where we can adopt Crandall's construction [C] of solutions to 
the Hamilton-Jacobi equation as test functions (see also [BJW]); (2) the comparison
principle of Hamilton-Jacobi equations without $u$-dependence (see [CIL] or [BJW]).

Since $f(x,p)=\sum_{i=1}^p p_i^2 (\ge 0)\in C^2(\Om\times R^m)$ satisfies both (f1) and (f2).
We have, as a consequence of theorem A,
\ss
\nind{\bf Corollary B}. {\it Suppose that $u\in W^{1,\fy}_{\hbox{cc}}(\Om)$ is an AMLE,
with respect to the Carnot-Carath\'edory metric $d_{\hbox{cc}}$. Then $u$ is a viscosity
solution to the subelliptic $\fy$-Laplacian equation (1.7).} 

Now we turn to the discussion on the uniqueness problem of absolute minimizers
of $F(\cdot, \Om)$ or viscosity solutions to eqn.(1.6). Although
the uniqueness might be true for general vector fields satisfying H\"ormander's
condition, we restrict our attention to the case where the vector fields generating
the Carnot-Carath\'edory metrics are horizontal vector fields associated
with a Carnot group $\bf G$. 

To describe the uniqueness results,  we recall that a Carnot group of step $r\ge 1$
is a simply connected Lie group $\bf G$ whose Lie algebra ${\sl g}$ admits a vector space
decomposition in $r$ layers ${\sl g}=V^1+V^2+\cdots +V^r$ having two properties:
(i) ${\sl g}$ is stratified, i.e., $[V^1, V^j]=V^{j+1}, j=1,\cdots, r-1$;
(ii) ${\sl g}$ is $r$-nilpotent, i.e. $[V^j, V^r]=0, j=1, \cdots, r$.  
$V^1$ is called the {\it horizontal} layer and $V^j, j=2,\cdots, r$, 
are {\it vertical} layers.
It is well-known (cf. Folland-Stein [FS]) 
that the exponential map, $\hbox{exp}: {\sl g} \to {\bf G}$, is a global
differmorphism so that we can identify $\bf G$ with ${\sl g}\equiv R^n$ via exp.
and $\bf G$ has an exponential coordinate system, here
$n=\hbox{dim}(\bf G)$ is the dimension of $\bf G$.
More precisely, Let $X_{i,j}$ for $1\le i\le m_j=\hbox{dim}(V^j)$ be a basis
of $V^j$ for $1\le j\le r$, which is orthonormal with respect to
an arbitrarily chosen Euclidean norm $\|\cdot\|$ on $\sl g$, with
respect to which the $V^j$'s are mutually orthogonal.
Then $p\in\bf G$ has coordinate $(p_{ij})_{1\le i\le m_j,
1\le j\le r}$ if $p=\hbox{exp}(\sum_{j=1}^r\sum_{i=1}^{m_j}(p_{ij}X_{i,j})$.
Let ${\cdot}$ denote the group multiplication on ${\bf G}$. Then it is known ([FS]) that
the group law $(x,y)\to x{\cdot} y$ is a polynomial map with respect to the exponential map.  
From now on, we set $m=m_1=\hbox{dim}(V^1)$ and denote
$X_i=X_{i,1}$ for $1\le i\le m$.
Two bi-Lipschitz equivalent metrics, on $\bf G$, we need are: (1) the Carnot-Carath\'edory 
metric $d_{\hbox{cc}}$ on $\bf G$
generated by $\{X_i\}_{i=1}^m$; (2) the gauge metric $d$ on $\bf G$ given as follows.
For $p=(p_{ij})_{1\le i\le m_j, 1\le j\le r}$, 
$$\|p\|_{\bf G}^{2r!}=\sum_{j=1}^r(\sum_{i=1}^{m_j}|p_{ij}|^2)^{r!\over j},$$
with the induced gauge distance
$$d(x,y)=\|x^{-1}y\|_{\bf G},  \ \ \forall x, y \in {\bf G}.$$
satisfying the invariant property
$$d(z\cdot x,z\cdot y)=d(x,y),  \ \ \forall x, y, z\in G.$$

Now we mention the Heisenberg group ${\bf H}^n$, which is
the simplest Carnot group of step two.  
${\bf H}^n\equiv {\bf C}^n\times R$ endowed with the group law:
for $(z_1,\cdots, z_n, t), (z_1',\cdots, z_n',t')\in {\bf C}^n\times R$
$$(z_1,\cdots,z_n, t)\cdot (z_1', \cdots, z_n', t')
=(z_1+z_1',\cdots, z_n+z_n', t+t'+2\hbox{Im}(\sum_{i=1}^n z_i\bar {z_i'})),$$
whose Lie algebra ${\sl h}=V_1+V_2$ with $V_1={\hbox{span}}\{X_i, Y_i\}_{1\le i \le n}$
and $V_2={\hbox{span}}\{T\}$, where
$$T = 4{\del\over\del t}, \ X_i={\del\over\del x_i}-2y_i{\del\over\del t},\ \ Y_i={\del\over\del 
y_i}+2x_i{\del\over\del t},
\ \ 1\le i\le n.$$

A function $f\in C^2(R^m)$ is 
strictly convex if there is a $C_0>0$ such that $D^2 f\ge C_0$.
Now we are ready to state the uniqueness theorem.   
\ss
\nind{\bf Theorem C}. {\it Let $\bf G$ be a Carnot group
and $\Om\ci G$ be a bounded domain. 
Assume that $f\in C^2(R^m, R_+)$ is strictly convex, homogeneous of
degree $\alpha\ge 1$, and $f(p)>0$ for $p\not=0$. 
Then, for any $\phi\in W^{1,\fy}_{\hbox{cc}}(\Om)$, the Dirichlet problem
$$\eqalignno{A(Xu, (D^2u)^*):=
-\sum_{i j=1}^m f_{p_i}(Xu)f_{p_j}(Xu)X_i X_j u &= 0,  \ \hbox{ in }\ \Om, &(1.13)\cr
u & = \phi,  \ \hbox{ on } \del\Om.\cr}$$
has at most one viscosity solution in $C(\bar\Om)$.}

Although the operator $A$ is degenerate subelliptic, one can check that the operator
${\bar A}(x, Du,D^2u)\equiv A(Xu, (D^2u)^*)$ has $x$-dependence and is not degenerate
elliptic (see [CIL] for its definition). 
Therefore, the uniqueness theorems, by Jensen [J1], Ishii [I], or Jensen
-Lions-Souganidis [JLS], on viscosity solutions to 2nd order elliptic PDEs, are not 
applicable directly here. Our ideas are: (i) We observe that eqn.(1.13) is invariant
under group multiplications: for any $a\in \bf G$, if $u\in C({\bf G})$ is
a viscosity solution to eqn.(1.13), then $u_a(x)=u(a\cdot x): {\bf G}\to R$
is also a viscosity solution to eqn.(1.13). This enables us to extend the
sup/inf convoluation construction by [JLS] to $\bf G$ to convert
viscosity sub/supersolutions of eqn.(1.13) into semiconvex/concave 
sub/supersolutions. (ii) We modify Jensen's original arguments
[J1] [J2] to prove a comparison principle between semiconvex
subsolutions and semiconcave {\it strict} supersolutions to any
2nd order degenerate subelliptic equations, which is valid for
any vector fields satisfying H\"ormander's condition. (iii) We adopt
Jensen's approximation scheme by $p$-Laplacians ([J2]) to build
viscosity solutions to two auxiiary equations, with horizontal gradient constraints,
having the properties that any supersolution can be converted into {\it strict}
supersolution under small perturbations. (iv) Finally, we apply the 
comparison principle for the two auxiliary equations to prove the uniqueness of eqn.(1.13).  

As a consequence of theorem A and C,  we have

\nind{\bf Corollary D}. {\it  Let $\bf G$ be a Carnot group and $\Om\ci \bf G$ be
a bounded domain. Assume that $f\in C^2(R^m, R_+)$ is strictly convex, homogeneous of
degree $\alpha\ge 1$, and $f(p)>0$ for $p\not=0$. 
Then, for any $\phi\in W^{1,\fy}_{\hbox{cc}}(\Om)$, there
is a unique absolute minimizer $u\in W^{1,\fy}_{\hbox{cc}}(\Om)$, with
$u|_{\del\Om}=\phi$, to the functional $F(v,\Om)=\sup_{x\in\Om}f(Xv)$,
and the eqn.(1.13) has a unique viscosity solution in $C(\bar\Om)$.
In particular, $\phi$ has a unique \hbox{AMLE} in $W^{1,\fy}_{\hbox{cc}}(\Om)$ and
the subelliptic $\fy$-Laplacian eqn. (1.7) has a unique viscosity solution.} 

We would like to remark that Bieske [B1] [B2] has previously proved the uniqueness
of both AMLE and viscosity solution to eqn. (1.7) for Heisenberg group ${\bf H}^n$
and Grushin type plane. However, our methods are considerably different. 
Manfredi, in a forthcoming paper [M],
studies some uniqueness issues for uniformly subelliptic 2nd order PDEs on Carnot groups.

The paper is written as follows. In \S2, we outline the proof of theorem A.
In \S3, we discuss the sup/inf convolution construction on Carnot group $\bf G$.
In \S4, we discuss the comparison principle between semiconvex subsolutions and
{\it strict} semiconcave supersolutions to any degenerate subelliptic equations
associated to vector fields satisfying H\"ormander's condition. In \S5, we 
study two auxiliary equations to eqn. (1.13), 
with horizontal gradient constraints. In \S6, we prove theorem C.

\bs
\nind {\S2}. Proof of theorem A
\ss
This section is devoted to the proof of theorem A. It contains two steps: (i) the construction
of test functions by solving the Hamilton-Jacobi equation, which is motivated by [BJW] and [C];
(ii) the comparison between viscosity subsolution and classical strict supersolution of
the Hamilton-Jacobi equation, which is motivated by [CIL] and [BJW].  

\ss
\nind{\bf Proof of theorem A}. 
It suffices to prove that if $u$ fails to be a viscosity 
subsolution of eqn.(1.12) at the point $x=0\in\Om$ then $u$ fails to
be an absolute minimizer of $F(\cdot,\Om)$. This assumption implies that
there is an $r_0>0$ and $\phi\in C^2(\Om)$ for which
$B_{r_0}(0)\ci\ci\Om$ such that
$$0=u(0)-\phi(0)\ge u(x)-\phi(x), \ \forall x\in \Om, \eqno(2.1)$$
but
$$-\sum_{i=1}^m X_i(f(x,X\phi))f_{p_i}(x,X\phi)(0)=C_0>0. \eqno(2.2)$$
Now we have
\ms
\nind{\bf Lemma 2.1}. {\it There exist a neighborhood
$V$ of $0$ and an $\Phi \in C^2(V)$ such that
$$ \Phi(0)=\phi(0),\ \ D\Phi(0)=D\phi(0), \ \ D^2\Phi(0)>D^2\phi(0), \eqno(2.3)$$
and
$$f(x, X\Phi(x))=f(0,X\phi(0))>0, \ \forall x\in V. \eqno(2.4)$$} 
\ss
\nind{\bf Proof}. Since $\{X_i\}_{i=1}^m\ci C^2(\Om)$,
there is $(a_{ij})_{1\le i\le m, 1\le j\le n}\in C^2(\Om, R^{mn})$ such that
$$X_i(x)=\sum_{j=1}^n a_{ij}(x){\del\over\del x_j}, \ \forall x\in\Om.$$ 
Define ${\bar f}:\Om\times R^n\to R$ by
$${\bar f}(x,q_1,\cdots, q_n)=f(x, \sum_{j=1}^n a_{1j}(x)q_j, \cdots, \sum_{j=1}^n a_{mj}(x)q_j),
 \ \ \forall (x, q_1,\cdots,q_n)\in\Om\times R^n.$$
Note that, for any $(x,q)\in \Om\times R^n$ and $1\le i\le n$, we have
$${\del {\bar f}\over\del q_i}(x,q)=\sum_{k=1}^m a_{ki}(x)
{\del f\over\del p_k}(x,\sum_{j=1}^n a_{1j}(x)q_j,\cdots, \sum_{j=1}^n a_{mj}(x)q_j).$$
Moreover, since $\phi\in C^2(\Om)$, it is easy to see that 
$${\bar f}(x, D\phi(x))=f(x,X_1\phi(x), \cdots, X_m\phi(x))=f(x, X\phi(x)),
\ \forall x\in \Om. \eqno(2.5)$$
Therefore, for any $x\in\Om$, we have
$$\sum_{j=1}^m X_j(f(x,X\phi(x)))f_{p_j}(x,X\phi(x))
=\sum_{i=1}^n {\del\over\del x_i}({\bar f}(x,D\phi(x)){\del{\bar f}\over\del q_i}(x,D\phi(x)).
\eqno (2.6)$$
This, combined with (2.2), implies
$$ A(0, D\phi(0), D^2\phi(0)):=-\sum_{i=1}^n {\del\over\del x_i}({\bar f}(x,D\phi))
{\del {\bar f}\over\del q_i}(x, D\phi)(0)=C_0 >0. \eqno(2.7)$$
Now we can apply exactly the step one of Crandall's argument ([C], page 275-276)
to conclude that there are a neighborhood $V$ of $0$ and an
$\Phi \in C^2(V)$ such that
$$\Phi(0)=\phi(0), \ \ D\Phi(0)=D\phi(0), \ \ D^2\Phi(0)>D^2\phi(0),$$
and
$${\bar f}(x, D\Phi(x))={\bar f}(0,D\phi(0)), \ \ \forall x\in V. \eqno(2.8)$$   
(2.8), combined with (2.5), gives (2.4). 

To see $f(0,X\phi(0))>0$, we observe that (2.2) implies
$$f_p(0,X\phi(0)):=({\del f\over\del p_1}(0,X\phi(0)), \cdots, {\del f\over\del p_m}(0,X\phi(0)))
\not=0.$$
This, combined with the fact that $f_p(0,0)=0$, 
implies $X\phi(0)\not=0$. Note that the homogenity of $f$ implies that
$f(0,0)=0$. Therefore, $f(0,X\phi(0))>0$.  
This finishes the proof of Lemma 2.1.  \qed 
\ss
It follows from Lemma 2.1 that there exists an open neighborhood $V_1\ci V$
of $0$ such that $\Phi(x)>\phi(x)\ge u(x)$ for any $0\not=x\in V_1$. Since
$\Phi(0)=\phi(0)=u(0)$. Therefore, for any small $\epsilon>0$, there
exists another neighborhood $V_\epsilon\ci V_1$ of $0$ such that 
$$\Phi(x)-\epsilon < u(x),\  \ \forall x\in V_\epsilon;\ \ 
 \Phi(x)-\epsilon =u(x), \ \ \forall
x\in\del V_\epsilon. \eqno(2.9)$$
It follows from the absolute minimality of $u$ to $F(\cdot, \ \Om)$ that
$$F(u, V_\epsilon)\le F(\Phi-\epsilon, V_\epsilon)
=\sup_{x\in V_\epsilon}f(x,X\Phi(x))=f(0,X\phi(0)). \eqno(2.10)$$

Now we want to show that $u$ is a viscosity subsolution of the Hamilton-Jacobi 
equation (2.8) on $V_\epsilon$. More precisely, we have
\ss
\nind{\bf Lemma 2.2}. {\it Under the same notations as above.
$u\in W^{1,\fy}_{\hbox{cc}}(V_\epsilon)$ is a viscosity subsolution
to the Hamilton-Jacobi equation
$$f(x, Xu(x))-f(0,X\phi(0))=0, \ \  \hbox{ in }\ V_\epsilon.  \eqno(2.11)$$}
\ss
\nind{\bf Proof}. For any subdomain $U\ci\ci V_\epsilon$
and  $0<\delta<\hbox{dist}(U,\del V_\epsilon)$, here dist
denotes the euclidean distance.
Let $g_\delta: U\to R$ be
the usual $\delta$-mollifier of $g$ for any function
$g$ on $V_\epsilon$. Since $u\in W^{1,\fy}_{\hbox{cc}}(V_\epsilon)$,
$u_\delta$ converges uniformly to $u$ on $U$ as $\delta\rightarrow 0$. 
Since $f$ is quasiconvex in its 2nd variable by (f1), 
it follows from the Jensen inequality for quasiconvex functions (cf. [BJW] theorem 1.1)
that for any $x\in U$
$$f(x, (Xu)_\delta(x))\le F(u,V_\epsilon)\le f(0, X\phi(0)). $$
Hence 
$$\sup_{x\in U}f(x, (Xu)_\delta(x))\le f(0, X\phi(0)). \eqno(2.12)$$
On the other hand, for any $1\le i\le m$ and $x\in U$, we can estimate
$(X_iu)_\delta(x)-X_i(u_\delta)(x)$ as follows
$$\eqalignno{&(X_iu)_\delta(x)-X_i(u_\delta)(x)\cr
&=\int_{R^n}\eta_\delta(x-y)(\sum_{j=1}^n a_{ij}(y){\del\over\del y_j})(u(y)-u(x))\,dy\cr
&-\int_{R^n}\sum_{j=1}^n a_{ij}(x){\del\eta_\delta(x-y)\over\del x_j} (u(y)-u(x))\,dy\cr
&=\sum_{j=1}^n \int_{R^n}\{-{\del\over\del y_j}(a_{ij}(y)\eta_\delta(x-y))
-a_{ij}(x){\del\eta_\delta(x-y)\over\del x_j}\} (u(y)-u(x))\,dy\cr
&=\sum_{j=1}^n\int_{R^n}(a_{ij}(y)-a_{ij}(x)){\del\eta_\delta(x-y)\over\del x_j} (u(y)-u(x))\,dy\cr
&+\sum_{j=1}^n \int_{R^n}{\del a_{ij}(y)\over\del y_j}\eta_\delta(x-y)(u(y)-u(x))\,dy.\cr}$$
Therefore we have
$$\eqalignno{|(X_i u)_\delta(x)-X_i(u_\delta)(x)|&\le C
\max_{1\le j \le n}\|Da_{ij}\|_{L^\fy(\Om)}\int_{R^n}\{|\eta_\delta(x-y)|u(y)-u(x)|\cr
&+|y-x||D\eta_\delta(x-y)||u(y)-u(x)|\}\,dy\cr
&\le C\|X_i\|_{C^1(\Om)}\sup_{\|y-x\|\le\delta}|u(y)-u(x)|\cr
&\le C\|X_i\|_{C^1(\Om)}\|u\|_{W^{1,\fy}_{\hbox{cc}(V_\epsilon)}}\delta^{1\over r},\cr}$$
where $r\ge 1$ is the step of H\"ormander's condition.
This implies
$$\eqalignno{f(x, X(u_\delta)(x))&\le \sup_{x\in U}f(x, (Xu)_\delta(x))
+\|f_p\|_{L^\fy}\|X(u_\delta)-(Xu)_\delta\|_{L^\fy(U)}.\cr
&\le f(0,X\phi(0))+ C\delta^{2\over r},\ \ \forall x\in U. &(2.13)\cr}$$
This, combined with the compactness theorem for viscosity solutions (cf. [CIL]), yields
that $u$ is a viscosity subsolution to the eqn. (2.8) in $U$. Since $U$ exhausts $V_\epsilon$
as $\delta\rightarrow 0$, we have that $u$ is a viscosity subsolution of the eqn. (2.8)
in $V_\epsilon$. The proof of Lemma 2.2 is complete.   \qed  
\ss

Now we continue the proof of theorem A. It follows from (f2) that
$$f(x,(1+t)p)=(1+t)^\alpha f(x, p)=(1+g(t))f(x,p), \ \forall t>0, 
\ \forall (x,p)\in \Om\times R^n,$$
where $g(t)\equiv(1+t)^\alpha-1>0$ for $t>0$, for $\alpha\ge 1$.
This, combined with (2.8), implies that, for any $t>0$,
$$f(x, X((1+t)\Phi_\epsilon)(x))=(1+g(t))
f(0,X\phi(0))=f(0, X\phi(0))+\delta(t), \ \forall x\in V_\epsilon, \eqno(2.14) $$
where $\Phi_\epsilon\equiv\Phi-\epsilon$ and $\delta(t)=g(t)f(0,X\phi(0))>0$.
Therefore, for any $t>0$,
$(1+t)\Phi_\epsilon$ is a {\it strict}, {\it classical} supersolution of eqn. (2.8).
We can now apply the comparison theorem for the Hamilton-Jacobi eqn.(2.8)
(see, e.g. Crandall-Ishii-Lions [CIL]) to conclude that
$$\sup_{V_\epsilon}(u-(1+t)\Phi_\epsilon)\le \sup_{\del V_\epsilon}(u-(1+t)\Phi_\epsilon),
\ \forall t>0.  \eqno(2.15)$$
Taking $t$ into zero, we  have
$$\sup_{V_\epsilon}(u-\Phi_\epsilon)\le \sup_{\del V_\epsilon}(u-\Phi_\epsilon)=0.$$
This implies
$$u(x)\le \Phi_\epsilon(x), \ \forall x\in V_\epsilon.$$ 
This clearly contradicts with (2.9). Therefore the proof of
theorem A is complete.  \qed

\bs
\nind \S3. The construction of sup/inf convolutions on $\bf G$ 
\ss
This section is devoted to the construction of sup/inf convolutions on the
Carnot group $\bf G$, which is the necessary extension of Jensen-Lions-Souganidis
[JLS] we need for the proof of theorem C.

Let $\Om\ci\bf G$ be a bounded domain and $d: {\bf G}\times {\bf G}\to R_+$ be
the gauge distance defined in \S1. For any $\epsilon>0$, define
$$\Om_\epsilon=\{x\in \Om: \inf_{y\in {\bf G}\setminus\Om} d^{2r!}(x^{-1},y^{-1})
\ge\epsilon\}.$$
\ss
\nind{\bf Definition 3.1}. For any $u\in C(\bar\Om)$ and $\epsilon>0$, the
sup involution, $u_\epsilon$,  of $u$ is defined by
$$u^\epsilon(x)=\sup_{y\in\bar\Om}(u(y)-{1\over 2\epsilon}d(x^{-1},y^{-1})^{2r!}), 
\ \forall x\in\Om. \eqno(3.1)$$
Similarly, the inf involution, $v_\epsilon$, of $v\in C(\bar\Om)$ is defined by
$$v_\epsilon(x)=\inf_{y\in\bar\Om}(v(y)+{1\over 2\epsilon}d(x^{-1},y^{-1})^{2r!}), 
\ \forall x\in\Om. \eqno(3.2)$$

For $x\in \bf G$, let $|x|:=(\sum_{j=1}^r\sum_{i=1}^{m_j}x_{ij}^2)^{1\over 2}$ 
denote its euclidean norm. We recall
\ss
\nind{\bf Definition 3.2}. A function $u\in C(\bar\Om)$ is called semiconvex, if
there is a constant $C>0$ such that $u(x)+C|x|^2$ is convex; $u$ is called
semiconcave if $-u$ is semiconvex. Note that, for $u\in C^2(\Om)$, if
$D^2u(x)\ge -C$ for $x\in\Om$ then $u$ is semiconvex, here $D^2u$ denotes
the (full) hessian of $u$. 

Now we have the generalized version of [JLS].
\ss
\nind{\bf Proposition 3.3}. {\it For $u, v\in C(\bar\Om)$, denote
$R_0=2\max\{\|u\|_{L^\fy(\Om)},\|v\|_{L^\fy(\Om)}\}$. Then,
for any $\epsilon>0$,  $u^\epsilon, v_\epsilon\in W^{1,\fy}_{\hbox{cc}}(\Om)$
satisfying:

\nind (1) $u^\epsilon$ is semiconvex and $v_\epsilon$ is semiconcave.

\nind (2) $u^\epsilon$ is monotonically nondecreasing w.r.t. $\epsilon$
and converges uniformly to $u$ on $\Om_{(1+4R_0)\epsilon}$;
$v_\epsilon$ is monotonically nonincreasing w.r.t. $\epsilon$ and 
converges uniformly to $v$ on $\Om_{(1+4R_0)\epsilon}$.

\nind (3) if $u$ (or $v$ respectively) is a viscosity subsolution (or
supersolution respectively) to a degenerate subelliptic
equation:
$$ B(Xu, (D^2u)^*) = 0 \ \hbox{ in }\ \Om,  \eqno(3.3)$$
then $u^\epsilon$ (or $v_\epsilon$) is a viscosity subsolution
(or supersolution respectively) to eqn. (3.3) in $\Om_{(1+4R_0)\epsilon}$.}
\ss
\nind{\bf Proof}. Since the proof of $v_\epsilon$ is similar
to that of $u^\epsilon$, we only prove the conclusions
for $u^\epsilon$. (1) Since $\Om\ci\bf G$ is bounded, it is easy to see 
from the formula of $d$ that
$$C_d(\Om)\equiv \|D^2_x(d(x^{-1},y^{-1})^{2r!})\|_{L^\fy(\bar\Om\times\bar\Om)}<\fy.$$
Therefore, for any $y\in\bar\Om$,
$$\tilde u^\epsilon_y(x):=u(y)-{1\over 2\epsilon}
d(x^{-1},y^{-1})^{2r!}+{C_d(\Om)\over 2\epsilon}|x|^2,
\ \forall x\in\Om,$$
has nonnegative hessian and is convex.
Since the maximum for a family of convex functions is still convex, this
implies that
$$u_\epsilon(x)+{C_d(\Om)\over 2\epsilon}|x|^2=\sup_{y\in\bar\Om}\tilde u^\epsilon_y(x)$$
is convex so that $u_\epsilon$ is semiconvex. It is well-known that semiconvex functions are
Lipschitz with respect to the euclidean metric so that $u^\epsilon\in W^{1,\fy}_{\hbox{cc}}(\Om)$.

(2) It is easy to see that for any $\epsilon_1<\epsilon_2$ 
$u^{\epsilon_1}(x)\le u^{\epsilon_2}(x)$ and $u(x)\le u^\epsilon(x)\le R_0$ for any $x\in\Om$. 
Observe that for any $x\in\Om$, $u^\epsilon(x)=\sup_{\bar\Om\cap\{d^{2r!}(x^{-1},y^{-1})
\le 4R_0\epsilon\}}
(u(y)-{1\over 2\epsilon}d(x^{-1},y^{-1})^{2r!})$. Therefore, for any $x\in 
\Om_{(1+4R_0)\epsilon}$,
$u^\epsilon(x)$ is attained at points $y\in \Om$. To see $u^\epsilon\rightarrow u$
uniformly on $\Om_{(1+4R_0)\epsilon}$, we observe that
if $u^\epsilon(x)$ is attained by $x_\epsilon$ then
$$u_{\epsilon\over 2}(x)\ge u(x_\epsilon)-{1\over\epsilon}d(x^{-1}, x_\epsilon^{-1})^{2r!}
=u_\epsilon(x)-{1\over 2\epsilon}d(x^{-1}, x_\epsilon^{-1})^{2r!}.$$
Hence
$$\lim_{\epsilon\rightarrow 0}{1\over\epsilon}d(x^{-1}, x_\epsilon^{-1})^{2r!}=0.$$
This implies that $x_\epsilon\rightarrow x$ and $\lim_{\epsilon\rightarrow 0}u^\epsilon(x)=u(x)$.
Moreover, since $$|u^\epsilon(x_1)-u^\epsilon(x_2)|\le |u(x_1)-u(x_2)|,\ \forall x_1, x_2\in\Om,$$
the convergence is uniform on $\Om_{(1+4R_0)\epsilon}$.

(3) For any $x_0\in\Om_{(1+4R_0)\epsilon}$, let $\phi\in C^2(\Om_{(1+4R_0)\epsilon})$ be such that
$$u_\epsilon(x_0)-\phi(x_0)\ge u_\epsilon(x)-\phi(x), \ \ \forall x\in\Om_{(1+4R_0)\epsilon}.$$
It follows from the proof of (2) above that there exists a $y_0\in\Om$ such that
$$u_\epsilon(x_0)=u(y_0)-{1\over 2\epsilon}d(x_0^{-1}, y_0^{-1})^{2r!}.$$
Therefore, we have
$$u(y_0)-{1\over 2\epsilon}d(x_0^{-1},y_0^{-1})^{2r!}-\phi(x_0)
\ge u(y)-{1\over 2\epsilon}d(x^{-1}, y^{-1})^{2r!}-\phi(x), \forall x, y\in\Om_{(1+4R_0)\epsilon}.$$
For $y$ near $y_0$, since $x=x_0\cdot y_0^{-1}\cdot y\in \Om_{(1+4R_0)\epsilon}$, we have
$$u(y_0)-\phi(x_0\cdot y_0^{-1}\cdot y_0)\ge u(y)-\phi(x_0\cdot y_0^{-1}\cdot y).$$
Set $\tilde{\phi}(y)=\phi(x_0\cdot y_0^{-1}\cdot y)$ for $y\in\Om_{(1+4R_0)\epsilon}$
near $y_0$. Then ${\tilde\phi}$ touches $u$ from above at $y=y_0$ and we have
$$B(X{\tilde \phi}, (D^2{\tilde \phi})^*)(y_0)\le 0. \eqno(3.4)$$
Now using the left-invariance of $X_i$, we know
$$X{\tilde \phi}(y)=(X\phi)(x_0\cdot y_0^{-1}\cdot y),\ \  (D^2(\tilde \phi))^*(y)
=(D^2 \phi)^*(x_0\cdot y_0^{-1}\cdot y).$$
This implies
$$B(X\phi(x_0), (D^2\phi)^*(x_0))\le 0.$$
Hence $u^\epsilon$ is a viscosity subsolution of eqn.(3.3) on $\Om_{(1+4R_0)\epsilon}$
and the proof 
of the proposition is complete.           \qed
\bs
\nind \S4. Comparison principle between semiconvex subsolutions and semiconcave supersolutions
\ss
In this section, we establish the comparison principle between
semiconvex subsolutions and semiconcave {\it strict }supersolutions
for any 2nd order subelliptic, possibly degenerate, PDE on 
the Carnot-Carath\'edory metric space generated by vector fields
satisfying H\"ormander's condition. The argument is inspired
by the well-known maximum principle for semiconvex functions,
due to Jensen [J1] [J2], on 2nd order elliptic PDEs. Here 
we assume that $\{X_i=\sum_{j=1}^n a_{ij}(x){\del\over\del x_j}\}_{i=1}^m$
is a set of vector fields
on $R^n$ satisfying H\"ormander's condition. The main proposition
of this section is 
\ss
\nind{\bf Proposition 4.1}. {\it For a bounded domain $\Om\ci  R^n$.
Suppose that $B\in C(\Om \times{\Cal S}^m)$ is
degenerate subelliptic. Assume that $u\in C(\bar\Om)$
is a semiconvex subsolution to
$$B(Xw, (D^2w)^*) = 0,  \ \hbox{ in  }\ \Om, \eqno(4.1)$$
and  $v\in C(\bar\Om)$ is a semiconcave supersolution to
$$B(Xw, (D^2w)^*) -\mu =0,  \ \hbox{ in } \ \Om, \eqno(4.2)$$
for some $\mu>0$. Then
$$\sup_{\Om}(u-v)\le\sup_{\del\Om}(u-v). \eqno(4.3)$$}
\ss
\nind{\bf Proof}. Suppose that (4.3) were false. Then
$$ \sup_{\Om}(u-v)>\sup_{\del\Om}(u-v),$$
so that $u-v$ achieves its maximum on $\bar\Om$ at a $x_0\in\Om$.
Since $u-v$ is semiconvex, it is well-known (cf. [J2] page 67) that
$$Du(x_0), Dv(x_0) \hbox{ both exist and are equal},$$
$$\eqalignno{u(x)-u(x_0)-\lan Du(x_0), x-x_0\ran &=O(|x-x_0|^2), &(4.4)\cr
 v(x)-v(x_0)-\lan Dv(x_0), x-x_0\ran &=O(|x-x_0|^2), &(4.5)\cr}$$
where $\lan\cdot, \cdot\ran$ and $|\cdot|$ denote the Euclidean inner product
and Euclidean norm.
Let $R_0=\hbox{dist}(x_0,\del\Om)>0$ be the euclidean distance from
$x_0$ to $\del\Om$ and $R_1>0$ be such that both (4.4) and (4.5) hold, with $|x-x_0|<R_1$.
Set $R_2=\min\{R_0,R_1\}>0$. Then, for any $\rho>0$, 
define the rescaled maps $u^\rho, v^\rho: B_{R_2\rho^{-1}}\to
R$ by
$$\eqalignno{u^\rho(x)&={1\over\rho^2}(u(x_0+\rho x)-u(x_0)-\rho\lan Du(x_0), x\ran),\cr
v^\rho(x)&={1\over\rho^2}(v(x_0+\rho x)-v(x_0)-\rho\lan Dv(x_0), x\ran),\cr}$$
where the Euclidean addition and scalar multiplication are used.
Then, it is easy to see 
$$0=(u^\rho-v^\rho)(0)\ge (u^\rho-v^\rho)(x), \ \ \ \forall x\in B_{R_2 \rho^{-1}}.$$
It follows from (4.4) and (4.5) that, for any $R>0$, 
there exists an $\rho_0=\rho_0(R)>0$ such that (i)
$\{u^\rho\}_{\{0<\rho\le \rho_0\}}$ are uniformly bounded,
uniformly semiconvex, and uniformly Lipschitz continuous in $B_R$;
(ii) $\{v^\rho\}_{\{0<\rho\le \rho_0\}}$ are uniformly bounded,
uniformly semiconcave, and uniformly Lipschitz continuous in $B_R$.
Therefore, by the Cauchy diagonal process, we may assume that
there is $\rho_i\downarrow 0$ such that $u^{\rho_i}\rightarrow u^*$,
$v^{\rho_i}\rightarrow v^*$ locally uniformly in $R^n$.
Moreover, it is not difficult to see that $u^*$ is locally bounded, semiconvex in $R^n$,
$v^*$ is locally bounded, semiconcave in $R^n$, and
$$0=(u^*-v^*)(0)\ge (u^*-v^*)(x), \ \ \ \forall x\in R^n.$$
Now we need 
\ss
\nind{\bf Claim 4.2}. {\it $u^*$ is a viscosity subsolution to
$$B_1(D^2 w)=0, \ \hbox{ in }\  R^n, \eqno(4.6)$$
and $v^*$ is a viscosity supersolution to
$$B_2(D^2 w)+\mu=0, \ \hbox{ in }\ R^n ,\eqno(4.7)$$
where $B_1, B_2: {\Cal S}^m\to R$ are defined by
$$B_1(M)=B(Xu(x_0),
\{\sum_{k,l=1}^n (a_{ik}(x_0)a_{jl}(x_0)M_{kl}
+a_{ik}(x_0){\del a_{jl}\over\del x_k}(x_0){\del u\over\del x_l}(x_0))\}_{1\le i,j \le m}), 
$$
$$B_2(M)=B(Xv(x_0), \{\sum_{k,l=1}^n (a_{ik}(x_0)a_{jl}(x_0)M_{kl}+
a_{ik}(x_0){\del a_{jl}\over\del x_k}(x_0){\del v\over\del x_l}(x_0)\}_{1\le i,j\le m})).$$ }

This claim follows from the compactness theorem (cf. [CIL]) for a family
of viscosity sub/supersolutions to 2nd order PDEs.
Since $u^\rho$ is a viscosity subsolution to  
$$B(Xu(x_0)+\rho X^\rho_w,
\{\sum_{k,l=1}^n a^\rho_{ik}a^\rho_{jl}{\del^2 w\over\del x_k\del x_l}
+a^\rho_{ik}({\del a_{jl}\over\del x_k})^\rho({\del u\over\del x_l}(x_0)
+\rho {\del w\over\del x_l})\}_{1\le i,j \le m})=0, \eqno(4.8)$$
and $v^\rho$ is a viscosity supersolution to
$$B(Xv(x_0)+\rho X^\rho_w,
\{\sum_{k,l=1}^n a^\rho_{ik}a^\rho_{jl}{\del^2 w\over\del x_k\del x_l}
+a^\rho_{ik}({\del a_{jl}\over\del x_k})^\rho({\del v\over\del x_l}(x_0)
+\rho {\del w\over\del x_l})\}_{1\le i,j \le m})=0, \eqno(4.9)$$
where $X^\rho=(X_1^\rho,\cdots, X_m^\rho)$, $X_i^\rho(x)=X_i(x_0+\rho x)$,
$a_{ik}^\rho(x)=a_{ik}(x_0+\rho x)$, and
$({\del a_{jl}\over\del x_k})^\rho(x)={\del a_{jl}\over\del x_k}(x_0+\rho x)$.
\ss
To see (4.8). Let $x_1\in B_{R_2\rho^{-1}}$ and $\phi\in C^2(B_{R_2\rho^{-1}})$ be such that
$$0=u^\rho(x_1)-\phi(x_1)\ge u^\rho(x)-\phi(x), \ \forall x\in B_{R_2\rho^{-1}}.$$
It is straightforward to see 
$$\phi^\rho(x)\equiv u(x_0)+\lan Du(x_0), x-x_0\ran +\rho^2\phi({x-x_0\over\rho}) $$ satisfies
$$0=u(x_0+\rho x_1)-\phi^\rho(x_0+\rho x_1)
\ge u(x)-\phi^\rho(x), \ \forall x\in B_{R_0}(x_0).$$
This, combined with the fact that $u$ is a viscosity subsolution to (4.6),
implies
$$B(X\phi^\rho, (D^2\phi^\rho)^*)(x_0+\rho x_1)\ge 0. \eqno(4.10)$$
Direct calculations yield
$${\del \phi^\rho\over\del x_k}(x)={\del u\over\del x_k}(x_0)
+\rho {\del\phi\over\del x_k}({x-x_0\over \rho}),$$
$${\del^2\phi^\rho\over\del x_k \del x_l}(x)
={\del^2\phi\over\del x_k \del x_l}({x-x_0\over \rho}).$$
Hence (4.10) implies (4.8).
It is clear that, after taking $\rho_i\rightarrow 0$, (4.8)-(4.9) imply (4.6)-(4.7).
This proves claim 4.2.

Since $u^*-v^*$ is semiconvex and achieves its maximum at $x=0$,
we can apply Jensen's maximum principle for semiconvex functions
(cf. [J1] [J2]) to conclude that there exists $x_*\in R^n$ such that
$D^2u^*(x_*), D^2v^*(x_*)$ both exist and $D^2(u^*-v^*)(x_*)\le 0$.
Let $M_1, M_2:{\Cal S}^m \to R$ be given by
$$M_1^{ij}=\sum_{k,l=1}^n (a_{ik}(x_0)a_{jl}(x_0){\del^2 u^*\over\del x_k\del x_l}(x_*)
+a_{ik}(x_0){\del a_{jl}\over\del x_k}(x_0){\del u\over\del x_l}(x_0)), 1\le i,j \le m, $$
and
$$M_2^{ij}=\sum_{k,l=1}^n (a_{ik}(x_0)a_{jl}(x_0){\del^2 v^*\over\del x_k\del x_l}(x_*)
+a_{ik}(x_0){\del a_{jl}\over\del x_k}(x_0){\del v\over\del x_l}(x_0)), 1\le i,j \le  m.$$
Since $Du(x_0)=Dv(x_0)$, we have, for any $p\in R^m$,
$$\sum_{1\le i, j\le m}(M_1^{ij}-M_2^{ij})p_ip_j =
\sum_{k,l=1}^n (\sum_{i=1}^m p_ia_{ik}(x_0))(\sum_{j=1}^m p_ja_{jl}(x_0))
{\del^2 (u-v)^*\over\del x_k\del x_l}(x_*)\le 0.$$
Hence $M_1\le M_2$. This, combined with the subellipticity of $B$ and $Xu(x_0)=Xv(x_0)$, implies
$$B(Xu(x_0), M_1)-B(Xv(x_0), M_2)\ge 0. \eqno(4.11)$$
On the other hand
$$B_1(D^2u^*(x_*))-B_2(D^2v^*(x_*))
=B(Xu(x_0), M_1)-B(Xv(x_0), M_2)\le -\mu<0. \eqno(4.12) $$
This contradicts with (4.11) and the proof of proposition 4.3 is complete.    \qed
\bs
\nind{\S5}. Auxiliary equations with horizontal gradient constraints
\ss
Due to the degenerancy of eqn.(1.13), we can't establish a comparison
principle for solutions to eqn.(1.13) directly. To get around the issue,
we follow Jensen's approximation scheme ([J2]) to  construct two auxiliary
equations with horizontal gradient constraints, to which supersolutions can
be deformed into {\it strict} supersolutions under small perturbations.
This section is valid for Carnot-Carath\'edory metric spaces associated
with vector fields satisfying H\"ormander's condition. In this section,
we assume that $\{X_i\}_{i=1}^m$ is a set of vector fields on $R^n$
satisfying H\"ormander's condition. First, we have
 
\ss
\nind{\bf Lemma 5.1}. {\it Suppose that $f\in C^2(R^m,R_+)$ is homogeneous of
degree $\alpha\ge 1$. Let $v\in C(\bar\Om)$
be a viscosity supersolution to
$$\min\{f(Xw)-\epsilon, \ \ -\sum_{ i j=1}^m f_{p_i}(Xw)f_{p_j}(Xw)X_iX_j w\}=0, 
\ \hbox{ in }\Om, \eqno(5.1)$$
where $\epsilon>0$. Then, for any $\delta>0$, 
there exist an $\mu=\mu(\alpha, \epsilon,\delta)>0$ and $v_\delta\in C(\bar\Om)$, 
with $\|v_\delta-v\|_{L^\fy(\Om)}\le\delta$, such that 
$v_\delta$ is a viscosity supersolution of 
$$\min\{f(Xw)-\epsilon, \ \  -\sum_{i j=1}^m f_{p_i}(Xw) f_{p_j}(Xw)X_i X_j w\} -\mu=0,
 \ \hbox{ in }\ \ \Om. \eqno(5.2)$$}
\ss
\nind{\bf Proof}. It is similar to that by Jensen [J2] 
(see also Juutinen [J] and Bieske [B1,2]). We sketch it here.
We look for $v_\delta=g_\delta(v)$, where $g_\delta\in C^\fy(R)$ is monotonically increasing such
that $g_\delta^{-1}\in C^\fy(R)$. To find $g_\delta$, let $x_0\in\Om$ and
$\phi\in C^2(\Om)$ touch $v_\delta$ from below at $x_0$. Let $\phi_\delta=g_\delta^{-1}(\phi)$. 
Then $\phi_\delta $ touches $v$ from below at $x_0$ and 
$$\min\{f(X\phi_\delta)-\epsilon, \ \  -\sum_{i j =1}^m f_{p_i}(X\phi_\delta)f_{p_j}(X\phi_\delta) 
X_iX_j\phi_\delta\}(x_0)\ge 0.$$
Since 
$$X_i\phi=g_\delta'(\phi_\delta)X_i\phi_\delta, X_iX_j\phi=g_\delta'(\phi_\delta)X_i X_j\phi_\delta
+g_\delta''(\phi_\delta)X_i\phi_\delta X_j\phi_\delta,$$
we have, by the $\alpha$-homogenity of $f$,
$$f(X\phi(x_0))
=f(g_\delta'(\phi_\delta)X\phi_\delta(x_0))=(g_\delta'(\phi_\delta))^\alpha f(X\phi_\delta(x_0))
\ge (g_\delta'(\phi_\delta))^\alpha\epsilon, \eqno(5.3)$$
and
$$\eqalignno{&-\sum_{i j=1}^mf_{p_i}(X\phi) f_{p_j}(X\phi) X_i X_j \phi(x_0)\cr
&=g_\delta'(\phi_\delta)^{3\alpha}(-\sum_{i j=1}^m f_{p_i}(X\phi_\delta) f_{p_i}(X\phi_\delta) 
X_i X_j\phi_\delta)(x_0)\cr
& -g_\delta'(\phi_\delta)^{2\alpha}g_\delta''(\phi_\delta)(\sum_{i=1}^m 
f_{p_i}(X\phi_\delta) X_i\phi_\delta)^2(x_0)\cr
&\ge - g_\delta'(\phi_\delta)^{2\alpha} g_\delta''(\phi_\delta)\alpha^2\epsilon^2, &(5.4)\cr}$$
provided that $g_\delta''(\phi_\delta)<0$, here we have used (5.3) and
the identity $\sum_{i=1}^m f_{p_i}(p)p_i=\alpha f(p)$ in the last step.
Let $C_0=4\|v\|_{L^\fy(\Om)}<\fy$ and define
$$g_\delta(t)=(1+\delta)t-{\delta\over 4C_0}t^2$$
for $|t|\le 2C_0$ and then extend this function suitably to a monotonically
increasing function on $R$. Since $g'(t)\ge 1+{\delta\over 2}$ and $g''(t)=-{\delta\over 2C_0}$
for $|t|\le C_0$, we have
$$f(X\phi)(x_0)\ge (1+{\delta\over 2})\epsilon,$$
and
$$-\sum_{i j=1}^m f_{p_i}(X\phi) f_{p_j}(X\phi) X_i X_j \phi(x_0)\ge {\delta 
\alpha^2\epsilon^2\over 2C_0}.$$
Therefore, if we choose $\mu=\min\{{\delta\epsilon\over 2}, 
{\delta\alpha^2\epsilon^2\over 2C_0}\}>0$, then 
$$\min\{f(X\phi)-\epsilon, \ \  -\sum_{i j=1}^m 
f_{p_i}(X\phi) f_{p_j}(X\phi) X_i X_j\phi\}(x_0)\ge \mu.$$
The proof of Lemma 5.1 is complete. 

Since the argument is similar, we state without proof the analogous
Lemma on viscosity subsolutions.
\ss
\nind{\bf Lemma 5.2}. {\it Suppose that $f\in C^2(R^m,R_+)$
is of homogeneous of degree $\alpha\ge 1$. Let $u\in C(\bar\Om)$ be a 
viscosity subsolution to
$$\max\{\epsilon-f(Xw), \ \  -\sum_{i j=1}^m f_{p_i}(Xw)f_{p_j}(Xw) X_iX_j w\}=0,
\ \hbox{ in }\ \Om,  \eqno(5.5)$$
where $\epsilon>0$. Then, for any $\delta>0$,
there are an $\mu=\mu(\alpha,\epsilon,\delta)>0$
and $u_\delta\in C(\bar\Om)$, with $\|u_\delta-u\|_{L^\fy(\Om)}\le\delta$,
such that $u_\delta$ is a viscosity subsolution to the equation
$$\max\{\epsilon-f(Xw),  \ \ -\sum_{i j=1}^m f_{p_i}(Xw)f_{p_j}(Xw) X_i X_j w \}=-\mu, 
\ \hbox{ in } \ \Om. \eqno(5.6)$$ }

We end this section with existences of viscosity solutions to eqn. (1.13), (5.3), (5.5).
For this, we need both convexity of $f$ and $f(p)>0$ for $p\not=0$. More precisely, 

\ss
\nind{\bf Theorem 5.3}. {\it Suppose that
$f\in C^2(R^m,R_+)$ is strictly convex, homogeneous of degree $\alpha\ge 1$,
and $f(p)>0$ for $p\not=0$. Then, for any $g\in 
W^{1,\fy}_{\hbox{cc}}(\Om)$, we have

\nind{(1)}. There exists a viscosity solution $u\in W^{1,\fy}_{\hbox{cc}}(\Om)$
to eqn.(1.13) such that $u|_{\del\Om}=g$.

\nind(2). There exists a viscosity solution $u_\epsilon\in W^{1,\fy}_{\hbox{cc}}(\Om)$ of 
eqn. (5.3) such that $u_\epsilon|_{\del\Om}=g$.

\nind(3). There exists a viscosity solution $v_\epsilon\in W^{1,\fy}_{\hbox{cc}}(\Om)$
of the eqn. (5.5) such that $v_\epsilon|_{\del\Om}=g$.

\nind(4). There exists a continuous, nondecreasing function $\beta: R_+\to R_+$,
with $\beta(0)=0$, such that 
$$\|u_\epsilon-v_\epsilon\|_{L^\fy(\Om)}\le \beta(\epsilon). \eqno(5.7)$$}
\ss
\nind{\bf Proof}. The proof is based on $L^k$ approximation,  which was first carried
out by [BDM], and then by Jensen [J2] for the $\fy$-Laplacian case (see also
[J] [B1, 2]). For completeness, we
outline it here. Since (1) follows from (2) with $\epsilon=0$ and (3)
can be done exactly in the way as (2), we only sketch (2)
and (4) as follows.
For $1<k<\fy$, let $u_p\in W^{1,k}_{\hbox{cc}}(\Om)$ be the unique minimizer to
the functional
$$F_k(v)=\int_\Om(f(Xv)^k-\epsilon^{k-1}v), \ \forall v\in W^{1,k}_{\hbox{cc}}(\Om), 
\hbox{ with } v|_{\del\Om}=g.$$
The existence of $u_k$ can be obtained by the direct method, due to both 
the convexity of $f$ and $\alpha$-homogeneity of $f$, i.e.
$f(p)=|p|^\alpha f({p\over|p|})\ge |p|^\alpha\min_{|z|=1}f(z)\ge C|p|^\alpha$.
It is easy to verify that $u_k$ satisfies the subelliptic $p$-Laplacian
equation
$$-\sum_{i=1}^mX_i^*(kf(Xu_k)^{k-1}f_{p_i}(Xu))=-\epsilon^{k-1}, \ \ \hbox{ in }\ \Om, 
\eqno(5.8)$$
in the sense of distributions, here $X_i^*$ is the adjoint of $X_i$. 
Let $Q$ denote the homogeneous dimension of $R^n$,
with respect to the vector fields $\{X_i\}_{i=1}^m$. Then it follows from the
Sobolev inequality (see, e.g., [HK]) that
$\{u_k\}_{k\ge Q+1}$ is bounded and equicontinuous. Therefore we may
assume, after taking possible subsequences, that there exist a
$u_\epsilon\in W^{1,\fy}_{\hbox{cc}}(\Om)$ such that
$$u_{k}\rightarrow u_\epsilon \ \hbox{ in }\ C^0(\bar\Om)\cap_{Q+1\le k<\fy }
W^{1,k}_{\hbox{cc}}(\Om).$$
It is easy to see that $u_\epsilon|_{\del\Om}=g$. To show that
$u_\epsilon$ is a viscosity solution to the eqn. (2.1)., we need
\ss
\nind{\bf Claim 5.4}. {\it For $k\ge Q+1$, $u_k\in C(\bar\Om)$ is a viscosity
solution to the eqn.(5.8).}
\ss
For simplicity, we only indicate that $u_k$ is a viscosity subsolution. 
For, otherwise, there are $x_0\in\Om$ and
$\phi\in C^2(\Om)$ such that
$$0=u_k(x_0)-\phi(x_0)\ge u_k(x)-\phi(x), \ \forall x\in\Om,$$
but
$$-\sum_{i=1}^m X_i^*(kf(X\phi)^{k-1}f_{p_i}(X\phi))(x_0)+\epsilon^{k-1}=-C_0<0. \eqno(5.9)$$
Then there exists an $\delta_0>0$ such that
$$-\sum_{i=1}^m X_i^*(kf(X\phi)^{k-1}f_{p_i}(X\phi))(x)+\epsilon^{k-1}\le -{C_0\over 2}<0, 
\ \ \forall x\in B_{\delta_0}(x_0).\eqno(5.10)$$
For any small $\delta>0$, there is a neighborhood $V_{\delta}(\ci B_{\delta_0}(x_0))$
of $x_0$ such that $\phi_\delta\equiv \phi-\delta$ satisfies
$$\phi_\delta(x)<u_k(x), \ \forall x\in V_\delta;  
\ \ \phi_\delta(x)=u_k(x), \ \forall x\in \del V_\delta.$$
Note that $\phi_\delta$ also satisfies (5.10).
Multiplying (5.8) by $u_k-\phi_\delta$ and integrating over $V_\delta$, we have
$$\sum_{i=1}^m\int_{V_\delta} kf(Xu_k)^{k-1}f_{p_i}(Xu_k) X_i (u_k-\phi_\delta)
=\epsilon^{k-1}\int_{V_\delta}(u_k-\phi_\delta).\eqno(5.11)$$
On the other hand, multiplying (5.10) by $(u_k-\phi_\delta)(\le 0)$ and integrating over
$V_{\delta}$, we have
$$\sum_{i=1}^m\int_{V_\delta} k(f(X\phi_\delta))^{k-1}f_{p_i}(X\phi_\delta)
X_i(u_k-\phi_\delta)>\epsilon^{k-1}\int_{V_\delta}(u_k-V_\delta). \eqno(5.12)$$
Subtracting (5.11) from (5.12), we obtain
$$0>
k\int_{V_\delta}\sum_{i=1}^m(f(Xu_k)^{k-1}f_{p_i}(Xu_p)
-f(X\phi_\delta)^{k-1}f_{p_i}(X\phi_\delta)) X_i(u_k-\phi_\delta))\ge 0,$$
this contradicts with the convexity of $f$.
This finishes the proof of Claim 5.4.

Now we show that $u_\epsilon$ is a viscosity subsolution to the eqn. (5.3).
Let $x\in\Om$ and $\phi\in C^2(\Om)$ be such that
$$0=u_\epsilon(x)-\phi(x)\ge u_\epsilon(y)-\phi(y), \ \forall y\in\Om.$$
We need to show
$$\min\{f(X\phi(x))-\epsilon,  \ -\sum_{i j=1}^m f_{p_i}(X\phi)
 f_{p_j}(X\phi) X_i X_j \phi(x)\}\le 0.$$
Since this is true if $f(X\phi(x))\le \epsilon$, we may assume that 
$f(X\phi(x))\ge (1+2\delta)
\epsilon$ for some $\delta>0$. 
We know that there exist $x_k\in\Om$ such that 
$(u_{k}-\phi)$ achieves its maximum at $x_k$ and $x_k\rightarrow x$. We
may also assume that, for $k$ sufficiently large,
$$f(X\phi(x_k))\ge (1+\delta)\epsilon.$$
It follows from claim 5.4 that
$$-\sum_{i=1}^m X_i^*(kf(X\phi)^{k-1}f_{p_i}(X\phi))(x_k)\ge -\epsilon^{k-1}.$$
After expansion 
and dividing both sides by $k(k-1)f(X\phi)^{k-2}(x_k)$, this gives
$$\eqalignno{\sum_{i j =1}^m f_{p_i} (X\phi)f_{p_j}(X\phi)  X_i X_j \phi(x_k) 
&\ge -{\epsilon\over k(k-1)}\{{\epsilon\over f(X\phi(x_k))}\}^{k-2}\cr
&+{f(X\phi(x_k))\over (k-1)}\sum_{i=1}^m X_i^*(f_{p_i}(X\phi))(x_k).\cr}$$
This, after taking $k$ into $\fy$, gives
$$\sum_{i j=1}^m f_{p_i}(X\phi) f_{p_j}(X\phi) X_iX_j\phi(x)\ge 0.$$
One can argue slightly differently that
$u_\epsilon$ is also a viscosity supersolution to the eqn. (5.3). This finishes the
proof of (2).

Since $v_\epsilon$ is a limit, as $k\rightarrow\fy$,
of the minimizers $v_k$ to
$$G_k(v)=\int_\Om f(Xv)^k+\epsilon^{k-1}v, \ \forall v\in W^{1,k}_{\hbox{cc}}(\Om),
\hbox{ with } u|_{\del\Om}=g,$$
$v_k$ satisfies
$$-\sum_{i=1}^mX_i^*(kf(Xv_k)^{k-1}f_{p_i}(Xv_k))=\epsilon^{k-1},
\ \hbox { in }\ \Om. \eqno(5.13)$$
Multiplying (5.11) and (5.13) by $(u_k-v_k)$, integrating over
$\Om$, and subtracting each other, we get
$$\eqalignno{&\int_\Om k(f(Xu_k)^{k-1}f_{p_i}(Xu_k)-f(Xv_k)^{k-1}
f_{p_i}(Xv_k))X_i(u_k-v_k)\cr
&\le 4\epsilon^{k-1}\|u_k-v_k\|_{L^1(\Om)}. &(5.14)\cr}$$
Now we need 
\ss
\nind{\bf Claim 5.5}. {\it If $f\in C^2(R^m)$ is strictly convex, then for any $p, q\in R^m$ 
$$(f^{k-1}(p)f_p(p)-f^{k-1}(q)f_p(q))\cdot (p-q)
\ge C|p-q|^{\alpha(k-1)+2}.\eqno(5.15)$$}

To see (5.15), we observe that
$$\eqalignno{&(f^{k-1}(p)f_p(p)-f^{k-1}(q)f_p(q))\cdot(p-q)\cr
&=k^{-1}\int_0^1 {d\over dt}(f^k)_p(tp+(1-t)q)\,dt \cdot (p-q)\cr
&\ge \sum_{i j=1}^m\int_0^1 f^{k-1}(tp+(1-t)q)f_{p_ip_j}(tp+(1-t)q)\,dt(p_i-q_i)(p_j-q_j) \cr
&\ge C^{k-1}\int_0^1 |tp+(1-t) q|^{\alpha(k-1)}\,dt |p-q|^2, &(5.16)\cr}$$
where we have used the strict convexity of $f$:
$$\sum_{i j=1}^m f_{p_ip_j}(v)p_i p_j\ge C_0|p-q|^2,  \ \forall p, q, v \in R^m,$$
the $\alpha$-homogenity of $f$ and the fact $f(p)>0$ for $p\not=0$:
$$f(v)=|v|^\alpha f({v\over |v|})\ge \min_{|z|=1}f(z) |v|^\alpha\ge C|v|^\alpha,
\ \forall v\in R^m,$$
for some $C>0$ depending only on $f$. Since
$$\int_0^1 |tp+(1-t)q|^{\alpha(k-1)}\,dt\ge C|p-q|^{\alpha(k-1)},$$
(5.16) implies (5.15).
Putting (5.15) into (5.14), we obtain
$$
kC^{k-1}\int_{\Om}|Xu_k-Xv_k|^{\alpha(k-1)+2}\le C\epsilon^{k-1}.$$
This, combined with the H\"older inequality, implies
$$\int_{\Om}|Xu_k-Xv_k|\le k^{-{1\over \alpha(k-1)+2}}(C\epsilon)^{k-1\over \alpha(k-1)+2}
|\Om|^{\alpha(k-1)+1\over \alpha(k-1)+2}.$$
Taking $k$ into $\infty$, we have
$$\|Xu_\epsilon -Xv_\epsilon\|_{L^1(\Om)}\le C\epsilon^{1\over \alpha}.\eqno(5.17)$$
In view of the fact that $u_\epsilon, v_\epsilon \in W^{1,\fy}_{\hbox{cc}}(\Om)$,
(5.17) together with the interpolation inequality and
the Sobolev inequality yield that the function $\beta$ must exist as asserted in (4).  
\qed

\bs
\nind \S6. Proof of theorem C
\ss
This section is devoted to the proof of theorem C. Henceforth we assume
that $\{X_i\}_{i=1}^m$ are horizontal vector fields in a bounded
domain $\Om$ of the Carnot group $\bf G$. Since we can identify
$\bf G$ with $R^n$, $n=\hbox{dim}(\bf G)$, via the exponential map,
the results in \S4 and \S5 are all applicable to $\bf G$. The idea
to prove theorem C is based on the sup/inf convolution and
the comparison principle for both equations (5.1) and (5.5).
\ss
\nind{\bf Lemma 6.1}. {\it Under the same assumptions as theorem C.
For any $\epsilon>0$, if $v\in C(\bar\Om)$
is a viscosity subsolution to the eqn.(5.1) and $w\in C(\bar\Om)$
is a viscosity supersolution to the eqn.(5.1), Then
$$\sup_{x\in\Om}(v-w)(x)=\sup_{x\in\del\Om}(v-w)(x). \eqno(6.1)$$}
\ss
\nind{\bf Proof}. Suppose that (6.1) were false. Then there is an $\delta_0>0$
such that
$$\sup_{x\in\Om}(v-w)(x)\ge \sup_{x\in\del\Om}(v-w)(x)+\delta_0.$$
For any $\delta\in (0, {\delta_0\over 2})$, it follows from Lemma 5.1 that
there are $w_\delta\in C(\bar\Om)$, with $\|w_\delta-w\|_{L^\fy(\Om)}\le\delta$,
and $\mu=\mu(\delta,\epsilon,\alpha)>0$, such that $w_\delta$ is a
viscosity supersolution to the eqn.(5.2). Moreover, we have
$$\sup_{x\in\Om}(v-w_\delta)\ge\sup_{x\in\del\Om}(v-w_\delta)(x)+{\delta_0\over 4}. \eqno(6.2)$$
Now we apply proposition 3.3 to conclude that for
any $\delta\in (0,{\delta_0\over 2})$ there are a semiconvex $v^\delta\in W^{1,\fy}_{\hbox{cc}}(\Om)$
and a semiconcave ${\tilde {w_\delta}}\in W^{1,\fy}_{\hbox{cc}}(\Om)$ such that
$$\lim_{\delta\rightarrow 0}\max\{\|v^\delta-v\|_{L^\fy(\Om_{C\delta})},
\ \|{\tilde {w_\delta}}-w_\delta\|_{L^\fy(\Om_{C\delta})}\}=0, \eqno(6.3)$$
where $\Om_{C\delta}$ is defined in \S3. Moreover, $v^\delta$ is a viscosity
subsolution to the eqn. (5.1) and ${\tilde {w_\delta}}$ is a viscosity supersolution
to the eqn. (5.2) on $\Om_{C\delta}$. Therefore, we can apply proposition 4.1
to conclude that
$$\sup_{\Om_{C\delta}}(v^\delta-{\tilde {w_\delta}})
=\sup_{\del\Om_{C\delta}}(v^\delta-{\tilde {w_\delta}}). \eqno(6.4)$$
Taking $\delta$ into zero, this yields
$$\eqalignno{\lim_{\delta\rightarrow 0}\sup_{\Om_{C\delta}}(v-w)&=
\lim_{\delta\rightarrow 0}\sup_{\Om_{C\delta}}[(v-v^\delta)+(v^\delta-{\tilde{w_\delta}})
+({\tilde{w_\delta}}-w_\delta)+(w_\delta-w)]\cr
&=\lim_{\delta\rightarrow 0}\sup_{\Om_{C\delta}}(v^\delta-{\tilde{w_\delta}})\cr
&=\lim_{\delta\rightarrow 0}\sup_{\del\Om_{C\delta}}(v^\delta-{\tilde{w_\delta}})\cr
&=\sup_{\del\Om}(v-w).\cr}$$
This yields the desired contradiction.  The proof is complete. \qed

Similarly, we have the comparison principle for the equation (5.5).
\ss
\nind{\bf Lemma 6.2}. {\it Under the same assumptions as theorem C.
For any $\epsilon>0$. Let $v\in C(\bar\Om)$ be
a viscosity subsolution to the eqn.(5.5) and $w\in C(\bar\Om)$ be a viscosity
supersolution to the eqn.(5.5). Then
$$\sup_{x\in\Om}(v-w)(x)=\sup_{x\in\del\Om}(v-w)(x).$$}
\ss
We are now in a position to prove a maximum principle for solutions of the 
eqn. (1.13) 
\ss
\nind{\bf Lemma 6.3}. {\it Under the same assumptions as theorem C. For
a given $\phi\in W^{1,\fy}_{\hbox{cc}}(\Om)$, assume that $v\in C(\bar\Om)$
is a viscosity subsolution to the eqn.(1.13) and $w\in C(\bar\Om)$ is a
viscosity supersolution to the eqn.(1.13) such that $v|_{\del\Om}=w|_{\del\Om}=\phi$.
Then
$$v(x)\le w(x), \ \forall x\in\Om. \eqno(6.5)$$}
\ss
\nind{\bf Proof}. Let $v^+$ be a viscosity solution to the eqn.(5.1) and
$w^-$ be a viscosity solution to the eqn.(5.5), with
$v^+|_{\del\Om}=w^-|_{\del\Om}=\phi$, obtained 
by theorem 5.3.  Since subsolutions of eqn.(1.13) are also subsolutions to eqn.(5.1)
and supersolutions to eqn.(1.13) are also supersolutions to eqn.(5.5), we can apply 
Lemma 6.1, 6.2 to conclude that
$$\sup_{\Om}(v-v^+)=\sup_{\del\Om}(v-v^+)=0, \ \sup_{\Om}(w^{-}-w)=\sup_{\del\Om}(w^{-1}-w)=0.$$
Hence we have
$$\sup_{\Om}(v-w)\le\sup_{\Om}(v^+-w^-)\le \beta(\epsilon),$$
where $\beta$ is given by theorem 5.3. Since $\epsilon$ is arbitrary, this implies
$$\sup_{\Om}(v-w)\le 0.$$
This finishes the proof of Lemma 6.3. \qed

It is obvious that Lemma 6.3 yields the conclusion of theorem C. Therefore, the proof of
theorem C is complete. \qed
\bs
\bs
\cl{\bf REFERENCES}
\ss
\nind{[A1]} G. Aronsson, {\it 
Extension of functions satisfying Lipschitz conditions}. Ark. Mat. 6 1967 551--561 (1967).

\nind{[A2]} G. Aronsson, {\it
On the partial differential equation 
$u\sb{x}{}\sp{2}\!u\sb{xx} +2u\sb{x}u\sb{y}u\sb{xy}+u\sb{y}{}\sp{2}\!u\sb{yy}=0$}. 
Ark. Mat. 7 1968 395--425 (1968). 

\nind{[B]} N. Barron, {\it Viscosity solutions and analysis in $L^\fy$}. Nonlinear
analysis, differential equations and control (Montreal, QC, 1998), Kluwer Acad. Publ. Dordrecht,
1999 (1-60).

\nind{[BDM]} T. Bhatthacharya, E. DiBenedetto, J. Manfredi, {\it Limits as
$p\rightarrow\fy$ of $\Delta_p u_p =f$ and Related Extremal Problems}.
Rend. Sem. Mat. Univ. Politec Torino, 1989, 15-68.

\nind{[B1]} T. Bieske, {\it On $\infty$-harmonic functions on the Heisenberg group}. 
Comm. Partial Differential Equations 27 (2002), no. 3-4, 727--761.

\nind{[B2]} T. Bieske, {\it Viscosity solutions on Grushin-type planes}. 
Illinois J. Math. 46 (2002), no. 3, 893--911.

\nind{[BC]} T. Bieske, L. Capogna, {\it The Aronsson-Euler equation for
absolutely minimizing Lipschitz extensions with repsect to Carnot-Carath\'eory
metrics}. Preprint.

\nind{[BJW]} N. Barron, R. Jensen, C. Y. Wang, {\it The Euler equation and absolute minimizers of 
$L\sp \infty$ functionals}. Arch. Ration. Mech. Anal. 157 (2001), no. 4, 255--283

\nind{[C]} M. Crandall, {\it An Efficient Derivation of the Aronsson Equation}.
Arch. Rational Mech. Anal. 167 (2003) 4, 271-279

\nind{[CE]} M. Crandall, L. C. Evans, {\it A remark on infinity harmonic functions}. 
Proceedings of the USA-Chile Workshop on Nonlinear Analysis, 123--129 (electronic), 
Electron. J. Differ. Equ. Conf., 6, Southwest Texas State Univ., San Marcos, TX, 2001. 

\nind{[CEG]} M. Crandall, L. C. Evans, R. Gariepy, {\it Optimal Lipschitz extensions and the 
infinity Laplacian}. Calc. Var. Partial Differential Equations 13 (2001), no. 2, 123--139. 

\nind{[CIL]} M. Crandall, H. Ishii, P. L. Lions, {\it User's guide to viscosity solutions
of second order partial differential equations}.
Bull. Amer. Math. Soc. (N.S.) 27 (1992), no. 1, 1--67.

\nind{[CL]} M. Crandall, P. L. Lions, {\it Viscosity solutions of Hamilton-Jacobi equations}. 
Trans. Amer. Math. Soc. 277 (1983), no. 1, 1--42. 

\nind{[E]} L. C. Evans, {\it Estimates for smooth absolutely minimizing Lipschitz extensions.}
Electron. J. Differential Equations 1993, No. 03, 9 pp. 

\nind{[FS]} G. Folland, E. Stein, Hardy spaces on homogeneous groups. Mathematical Notes, 28. 
Princeton University Press, Princeton, N.J., 1982. 

\nind{[FSS]} B. Franchi, R. Serapioni, F. Serra Cassano, {\it
Meyers-Serrin type theorems and relaxation of variational integrals
depending on vector fields}. Houston J. Math. 22 (1996), no. 4, 859--890.

\nind{[GN]} N. Garofalo, D. Nhieu, {\it Lipschitz continuity,
global smooth approximations and extension theorems for Sobolev functions
in Carnot-Carathodory spaces}. J. Anal. Math. 74 (1998), 67--97.

\nind{[HK]} P. Hajlasz, P. Koskela, {\it Sobolev met Poincar\'e}.
Mem. Amer. Math. Soc. 145 (2000), no. 688,

\nind{[I]} H. Ishii, {\it On existence and uniqueness of viscosity 
solutions of fully nonlinear second-order elliptic PDEs}. Comm. Pure Appl. Math.
42 (1989) 14-45.

\nind{[J]} P. Juutinen, {\it Minimization problems for Lipschitz functions via viscosity solutions}. 
Ann. Acad. Sci. Fenn. Math. Diss. No. 115 (1998),

\nind{[J1]} R. Jensen,  {\it The maximum principle for viscosity solutions of fully nonlinear
second order partial differential equations}. Arch. Rational Mech. Anal. 101 (1988), no. 1, 1--27. 

\nind{[J2]} R. Jensen, {\it Uniqueness of Lipschitz extensions: minimizing
the sup norm of the gradient}. Arch. Rational Mech. Anal. 123 (1993), no. 1, 51--74.   

\nind{[JLS]} R. Jensen,  P. L. Lions, P. Souganidis, {\it A uniqueness result for viscosity
solutions of second order fully nonlinear partial differential equations}. 
Proc. Amer. Math. Soc. 102 (1988), no. 4, 975--978. 

\nind{[LM]} P. Lindqvist, J. Manfredi, {\it The Harnack inequality for $\infty$-harmonic functions.}
Electron. J. Differential Equations 1995, No. 04, 5 pp. 

\nind{[M]} J. Manfredi, {\it Fully nonlinear subelliptci equations}. In preparation.

\nind{[NSW]} A. Nagel, E. Stein, S. Wainger, {\it  Balls and metrics defined by vector fields.
I. Basic properties}. Acta Math. 155 (1985), no. 1-2, 103--147.

\end